\documentclass[12pt]{article}
\usepackage[top=1in,bottom=1in,left=1in,right=1in]{geometry}
\usepackage{indentfirst}
\usepackage{amsfonts,amsmath,amsthm,amssymb,hyperref,enumerate,bm}
\usepackage{longtable}
\usepackage[caption=false]{subfig}
\usepackage{leftidx}
\usepackage{lineno}
\usepackage{color,xcolor}
\usepackage{array}
\usepackage{latexsym}
\usepackage{float}
\usepackage{color}
\usepackage{longtable,supertabular,booktabs}
\usepackage{rotating}
\usepackage{cases}
\usepackage{multicol,multirow}
\usepackage{epsfig}
\usepackage{cite}
\usepackage{fancyhdr}       % fuer Deckblatt
\usepackage{dsfont}
\usepackage{hyperref}
\hypersetup{
	colorlinks=true,
	linkcolor=red,
	filecolor=blue,
	urlcolor=red,
	citecolor=blue,
}

\newtheorem{theorem}{Theorem}[section]

\newtheorem{remark}{Remark}[section]

\newtheorem{example}{Example}[section]

\newtheorem{problem}{Problem}[section]
\newtheorem{question}{Question}[section]

\newtheorem{fact}{Fact}[section]

\def\G1{G^\mathcal{C}}

\newcommand{\rad}{\text{\rm rad}}
\newcommand{\Mod}{\text{\rm mod~}}

 \newenvironment{prof}{\trivlist
      \item[\hskip\labelsep
      {\itshape Proof.}]\normalfont}
      {\hspace*{\fill}$\Box$\endtrivlist}

\begin{document}
\title{On the representation of rational numbers via Euler's totient function}
\author{
Weilin Zhang$^{1,}$\footnote{E-mail\,$:$ weilin@m.scnu.edu.cn.}\quad
Fengyuan Chen$^{1,}$\footnote{E-mail\,$:$ 332248560@qq.com.}\quad
Hongjian Li$^{2,}$\footnote{Corresponding author. E-mail\,$:$ lhj@gdufs.edu.cn. Supported by the Project of Guangdong University of Foreign Studies (Grant No. 2024RC063).}\quad
Pingzhi Yuan$^{1,}$\footnote{E-mail\,$:$ yuanpz@scnu.edu.cn. Supported by the National Natural Science Foundation of China (Grant No. 12171163) and the Basic and Applied Basic Research Foundation of Guangdong Province (Grant No. 2024A1515010589).}\\
{\small\it  $^{1}$School of Mathematical Sciences, South China Normal University,}\\
{\small\it Guangzhou 510631, Guangdong, P. R. China} \\
{\small\it  $^{2}$School of Mathematics and Statistics, Guangdong University of Foreign Studies,}\\
{\small\it Guangzhou 510006, Guangdong, P. R. China} \\
}

\date{}
 \maketitle
\date{}

\noindent{\bf Abstract}\quad
Let $b>1$ be an odd positive integer and $k, l \in \mathbb{N}$. In this paper, we show that every positive rational number can be written as $\varphi(m^{2})/(\varphi(n^{2}))^{b}$ and $\varphi(k(m^{2}-1))/\varphi(ln^{2})$, where $m, n\in \mathbb{N}$ and $\varphi$ is the Euler's totient function. At the end, some further results are discussed.


\medskip
\noindent{\bf MR(2020) Subject Classification} 11A25, 11D85

\section{Introduction} \label{Sec1}

Let $\mathbb{N}$ be the set of all positive integers. For any $n\in \mathbb{N}$, the Euler's totient function $\varphi(n)$ is defined as the number of positive integers up to $n$ that are relatively prime to $n$. Suppose that
\begin{equation*}
	n = \prod_{i=1}^{s}p_{i}^{\alpha_{i}}, ~\text{where}~ p_{1}<p_{2}<\dots<p_{s} ~\text{are primes and}~ \alpha_{i}\in\mathbb{N},
\end{equation*}
is the standard factorization of $n$. It is well-known that (see, e.g., \cite[p. 20]{IR1990})
\begin{equation*}
	\varphi(n) = \varphi\left(\prod_{i=1}^{s}p_{i}^{\alpha_{i}}\right) = \prod_{i=1}^{s}p_{i}^{\alpha_{i}-1}(p_{i}-1).
\end{equation*}

Throughout this paper, we always suppose that $\rad(n) = \prod_{i=1}^{s}p_{i}$, and let $v_{p_{i}}(n) = \alpha_{i}$ denote the $p_{i}$-valuation of $n$. For $q = m/n\in\mathbb{Q^{+}}$, where $m, n\in\mathbb{N}$, the notation $v_{p_{i}}(q)$ means that $v_{p_{i}}(q) = v_{p_{i}}(m) - v_{p_{i}}(n)$. As an easy exercise of the above formula, we have
\begin{equation*}
	\dfrac{\varphi(m)}{\varphi(n)} = \dfrac{m}{n}
\end{equation*}
when $\rad(n) = \rad(m)$.

In 2020, Krachun and Sun \cite{KS2020} proved that every positive rational number can be written in the form $\varphi(m^{2})/\varphi(n^{2})$, where $m, n\in \mathbb{N}$. Li, Yuan and Bai \cite{LYB2021} proved that such representation is unique with some natural restrictions on $\gcd(m, n)$. Recently, Krachun and Sun's result has been generalized in the following two forms
\begin{equation*}
	\dfrac{\varphi(km^{r})}{\varphi(ln^{s})} ~~\text{and}~~ \dfrac{(\varphi(m^{r}))^{a}}{(\varphi(n^{s}))^{b}}, ~~\text{where}~~ k, l, r, s, a, b \in \mathbb{N}.
\end{equation*}

Of the two, the former has been clearly studied by Li, Yuan and Bai \cite{LYB2020}. Let $k, l, r, s\in \mathbb{N}$ with $\max\{r, s\} \geqslant 2$. They proved that every positive rational number can be written in the form $\varphi(km^{r})/\varphi(ln^{s})$, where $m, n \in \mathbb{N}$, if and only if $\gcd(r, s) = 1$ or $(k, l, r, s) = (1, 1, 2, 2)$. Moreover, if $\gcd(r, s) > 1$, then the proper representation of such representation is unique.

In 2022, T. T. Vu \cite{V2022} modified Krachun and Sun's proof and showed that if $\gcd(ar, bs) = 1$, then every positive rational number can be written in the form $(\varphi(m^{r}))^{a}/(\varphi(n^{s}))^{b}$, where $m, n \in\mathbb{N}$. In fact, Vu's proof shows that there exist infinitely many such positive integer pairs $(m, n)$. At the end of \cite{V2022}, Vu proposed the following open problem.

\begin{problem}
	Besides $(a, b, r, s) = (1, 1, 2, 2)$ and $(a, b, r, s)$ with $\gcd(ar, bs) = 1$, are there any other positive integer quadruples $(a, b, r, s)$ such that every positive rational number can be written in the form $(\varphi(m^{r}))^{a}/(\varphi(n^{s}))^{b}$, where $m, n\in\mathbb{N}$?
\end{problem}

In this paper, we answer the above problem by establishing the following theorem.

\begin{theorem} \label{Thm1}
	Let $b>1$ be an odd integer. Then every positive rational number can be written as $\varphi(m^{2})/(\varphi(n^{2}))^{b}$, where $m, n\in\mathbb{N}$.
\end{theorem}

Note that $q\in\mathbb{Q^{+}}$ if and only if $1/q\in\mathbb{Q^{+}}$. Let $a>1$ and $b>1$ be odd integers. Then Theorem \ref{Thm1} shows that $(a, 1, 2, 2)$ and $(1, b, 2, 2)$ are other positive integer quadruples that satisfy the conditions in the above problem. Further discussion of the above problem will occur in Section \ref{Sec3}.

Another main purpose of this paper is to prove the following result.

\begin{theorem} \label{Thm2}
	Let $k$ and $l$ be positive integers. Then every positive rational number can be written as $\varphi(k(m^{2}-1))/\varphi(ln^{2})$, where $m, n\in\mathbb{N}$.
\end{theorem}

Let $(k, l, t)$ be a positive integer triple. If we generalize Krachun and Sun's result in the form $\varphi(k(m^{2} - t))/\varphi(ln^{2})$, where $m, n\in\mathbb{N}$, then Krachun and Sun's result and Theorem \ref{Thm2} show that this form can express all positive rational number for the triples $(1, 1, 0)$ and $(k, l, 1)$. For a further study of this form, we will mention it at the end of this paper.

\section{Proofs and Examples}

In this section, we give our proofs in more detail and present some examples for Theorems \ref{Thm1}--\ref{Thm2}.

\emph{Proof of Theorem \ref{Thm1}.} Let $2 = p_{1} < p_{2} < p_{3} < \dots$ be the sequence of all primes. Let $Q_{0} = \{1\}$ and for every $l\in\mathbb{N}$, let $Q_{l} = \left\{\prod_{i=1}^{l}p_{i}^{\alpha_{i}}: \alpha_{i}\in\mathbb{Z}, 1\leqslant i \leqslant l\right\}$ be the set of all positive rational numbers that can be expressed as the product of integer powers of the first $l$ primes. Let $q$ be a positive rational number. If $q = 1$, then the statement is true for $m = n = 1$. For the case $q \neq 1$, we assume that $q\in Q_{s}\setminus Q_{s-1}$ for some $s\in\mathbb{N}$. We now give a construction for $m$ and $n$ as follows.
	
	\textbf{Step $1$:} By Dirichlet's theorem on primes in an arithmetic progression, we can choose a prime $q_{1}$ such that $q_{1} \equiv 1 \left(\Mod \left(p_{s}!\right)^{t}\right)$, where $t\in\mathbb{N}$ such that
	\begin{equation*}
		t(b-1) > \max_{1\leqslant j \leqslant s}\left\{v_{p_{j}}\left(\dfrac{1}{q}\prod_{i=1}^{s}(p_{i}-1)\right)\right\}.
	\end{equation*}
	Since $b$ is an odd integer, there are positive integers $x_{1}$ and $y_{1}$ such that
	\begin{equation*}
		(2x_{1} - 1) - (2y_{1}-1)b = \alpha_{1}(b-1), ~\text{where}~ \alpha_{1} = 0.
	\end{equation*}
	Let $m_{1} = q_{1}^{x_{1}}$ and $n_{1} = q_{1}^{y_{1}}$. Then we have
	\begin{equation*}
		\dfrac{\varphi(m_{1}^{2})}{(\varphi(n_{1}^{2}))^{b}} = \dfrac{\varphi(q_{1}^{2x_{1}})}{(\varphi(q_{1}^{2y_{1}}))^{b}} = \dfrac{q_{1}^{2x_{1}-1}(q_{1}-1)}{q_{1}^{(2y_{1}-1)b}(q_{1}-1)^{b}} = \dfrac{q_{1}^{\alpha_{1}(b-1)}}{(q_{1} - 1)^{b-1}}.
	\end{equation*}
	
	\textbf{Step $2$:} Let $q_{2}$ be the maximal prime factor of $(q_{1} - 1)/q_{1}^{\alpha_{1}}$. Since $b$ is an odd integer, there are positive integers $x_{2}$ and $y_{2}$ such that
	\begin{equation*}
		(2x_{2} - 1) - (2y_{2}-1)b = \alpha_{2}(b-1), ~\text{where}~ \alpha_{2} = v_{q_{2}}(q_{1} - 1).
	\end{equation*}
	Let $m_{2} = q_{2}^{x_{2}}$ and $n_{2} = q_{2}^{y_{2}}$. Then we have
	\begin{equation*}
		\dfrac{\varphi(m_{2}^{2})}{(\varphi(n_{2}^{2}))^{b}} = \dfrac{\varphi(q_{2}^{2x_{2}})}{(\varphi(q_{2}^{2y_{2}}))^{b}} = \dfrac{q_{2}^{2x_{2}-1}(q_{2}-1)}{q_{2}^{(2y_{2}-1)b}(q_{2}-1)^{b}} = \dfrac{q_{2}^{\alpha_{2}(b-1)}}{(q_{2} - 1)^{b-1}}.
	\end{equation*}
	
	\textbf{Step $3$:} Let $q_{3}$ be the maximal prime factor of $(q_{1} - 1)(q_{2} - 1)/(q_{1}^{\alpha_{1}}q_{2}^{\alpha_{2}})$. Since $b$ is an odd integer, there are positive integers $x_{3}$ and $y_{3}$ such that
	\begin{equation*}
		(2x_{3} - 1) - (2y_{3}-1)b = \alpha_{3}(b-1), ~\text{where}~ \alpha_{3} = v_{q_{3}}\left((q_{1} - 1)(q_{2} - 1)\right).
	\end{equation*}
	Let $m_{3} = q_{3}^{x_{3}}$ and $n_{3} = q_{3}^{y_{3}}$. Then we have
	\begin{equation*}
		\dfrac{\varphi(m_{3}^{2})}{(\varphi(n_{3}^{2}))^{b}} = \dfrac{\varphi(q_{3}^{2x_{3}})}{(\varphi(q_{3}^{2y_{3}}))^{b}} = \dfrac{q_{3}^{2x_{3}-1}(q_{3}-1)}{q_{3}^{(2y_{3}-1)b}(q_{3}-1)^{b}} = \dfrac{q_{3}^{\alpha_{3}(b-1)}}{(q_{3} - 1)^{b-1}}.
	\end{equation*}
	
	Continuing this procedure, we assume that $k$ is the maximal positive integer such that $q_{k} > p_{s}$. If $k>1$, then we will proceed to the following step $k$.
	
	\textbf{Step $k$:} Let $q_{k}$ be the maximal prime factor of $\prod_{i=1}^{k-1}(q_{i} - 1)/\prod_{i=1}^{k-1}q_{i}^{\alpha_{i}}$. Since $b$ is an odd integer, there are positive integers $x_{k}$ and $y_{k}$ such that
	\begin{equation*}
		(2x_{k} - 1) - (2y_{k}-1)b = \alpha_{k}(b-1), ~\text{where}~ \alpha_{k} = v_{q_{k}}\left(\prod_{i=1}^{k-1}(q_{i} - 1)\right).
	\end{equation*}
	Let $m_{k} = q_{k}^{x_{k}}$ and $n_{k} = q_{k}^{y_{k}}$. Then we have
	\begin{equation*}
		\dfrac{\varphi(m_{k}^{2})}{(\varphi(n_{k}^{2}))^{b}} = \dfrac{\varphi(q_{k}^{2x_{k}})}{(\varphi(q_{k}^{2y_{k}}))^{b}} = \dfrac{q_{k}^{2x_{k}-1}(q_{k}-1)}{q_{k}^{(2y_{k}-1)b}(q_{k}-1)^{b}} = \dfrac{q_{k}^{\alpha_{k}(b-1)}}{(q_{k} - 1)^{b-1}}.
	\end{equation*}
	
	Therefore, if $k>1$, then we have
	\begin{equation*}
		\dfrac{\varphi\left(\left(\prod_{i=1}^{k}m_{i}\right)^{2}\right)}{\left(\varphi\left(\left(\prod_{i=1}^{k}n_{i}\right)^{2}\right)\right)^{b}} = \prod_{i=1}^{k}\dfrac{\varphi(m_{i}^{2})}{(\varphi(n_{i}^{2}))^{b}} = \dfrac{\prod_{i=1}^{k}q_{i}^{\alpha_{i}(b-1)}}{\prod_{i=1}^{k}(q_{i}-1)^{b-1}} = \dfrac{1}{A_{0}},
	\end{equation*}
	where $A_{0}\in Q_{s}$ is a positive integer such that  $v_{p_{i}}(A_{0}) \geqslant t(b-1)$ for every $i = 1, 2, \dots, s$. Note that this conclusion also holds for $k = 1$.
	
	Next, we will continue with the following $s$ steps.
	
	\textbf{Step $k + 1$:} Let $q_{k+1} = p_{s}$ and $\alpha_{k+1} = v_{q_{k+1}}(qA_{0})$. Since
	\begin{equation*}
		v_{q_{k+1}}(A_{0})\geqslant t(b-1) > \max_{1\leqslant j \leqslant s}\left\{v_{p_{j}}\left(\dfrac{1}{q}\prod_{i=1}^{s}(p_{i}-1)\right)\right\} > -v_{q_{k+1}}(q),
	\end{equation*}
	it follows that
	\begin{equation*}
		\alpha_{k+1} = v_{q_{k+1}}(qA_{0}) =  v_{q_{k+1}}(q) + v_{q_{k+1}}(A_{0}) > 0.
	\end{equation*}
	Note that $b$ is an odd integer. If $2 \mid \alpha_{k+1}$, then there are positive integers $x_{k+1}$ and $y_{k+1}$ such that
	\begin{equation*}
		(2x_{k+1}-1) - (2y_{k+1}-1)b = \alpha_{k+1}.
	\end{equation*}
	Let $m_{k+1} = q_{k+1}^{x_{k+1}}$ and $n_{k+1} = q_{k+1}^{y_{k+1}}$. Then we have
	\begin{equation*}
		\dfrac{\varphi\left(m_{k+1}^{2}\right)}{\left(\varphi\left(n_{k+1}^{2}\right)\right)^{b}} = \dfrac{\varphi\left(q_{k+1}^{2x_{k+1}}\right)}{\left(\varphi\left(q_{k+1}^{2y_{k+1}}\right)\right)^{b}} = \dfrac{q_{k+1}^{2x_{k+1}-1}(q_{k+1} - 1)}{q_{k+1}^{(2y_{k+1}-1)b}(q_{k+1}-1)^{b}} = \dfrac{q_{k+1}^{\alpha_{k+1}}}{A_{1}},
	\end{equation*}
	where $A_{1} = (q_{k+1}-1)^{b-1}$. If $2 \nmid \alpha_{k+1}$, then $x_{k+1} = (\alpha_{k+1} + 1)/2$ is a positive integer. Let $m_{k+1} = q_{k+1}^{x_{k+1}}$ and $n_{k+1} = 1$. Then we have
	\begin{equation*}
		\dfrac{\varphi\left(m_{k+1}^{2}\right)}{\left(\varphi\left(n_{k+1}^{2}\right)\right)^{b}} = \dfrac{\varphi\left(q_{k+1}^{2x_{k+1}}\right)}{\left(\varphi\left(1^{2}\right)\right)^{b}} = q_{k+1}^{2x_{k+1}-1}(q_{k+1} - 1) = \dfrac{q_{k+1}^{\alpha_{k+1}}}{A_{1}},
	\end{equation*}
	where $A_{1} = 1/(q_{k+1}-1)$.
	
	\textbf{Step $k + 2$:} Let $q_{k+2} = p_{s-1}$ and $\alpha_{k+2} = v_{q_{k+2}}(qA_{0}A_{1})$. Since
	\begin{equation*}
		v_{q_{k+2}}(A_{0})\geqslant t(b-1) > \max_{1\leqslant j \leqslant s}\left\{v_{p_{j}}\left(\dfrac{1}{q}\prod_{i=1}^{s}(p_{i}-1)\right)\right\} > v_{q_{k+2}}(q_{k+1}-1)-v_{q_{k+2}}(q),
	\end{equation*}
	it follows that
	\begin{equation*}
		\alpha_{k+2} = v_{q_{k+2}}(qA_{0}A_{1}) \geqslant v_{q_{k+2}}(q) + v_{q_{k+2}}(A_{0}) - v_{q_{k+2}}(q_{k+1} - 1) > 0.
	\end{equation*}
	Note that $b$ is an odd integer. If $2 \mid \alpha_{k+2}$, then there are positive integers $x_{k+2}$ and $y_{k+2}$ such that
	\begin{equation*}
		(2x_{k+2}-1) - (2y_{k+2}-1)b = \alpha_{k+2}.
	\end{equation*}
	Let $m_{k+2} = q_{k+2}^{x_{k+2}}$ and $n_{k+2} = q_{k+2}^{y_{k+2}}$. Then we have
	\begin{equation*}
		\dfrac{\varphi\left(m_{k+2}^{2}\right)}{\left(\varphi\left(n_{k+2}^{2}\right)\right)^{b}} = \dfrac{\varphi\left(q_{k+2}^{2x_{k+2}}\right)}{\left(\varphi\left(q_{k+2}^{2y_{k+2}}\right)\right)^{b}} = \dfrac{q_{k+2}^{2x_{k+2}-1}(q_{k+2} - 1)}{q_{k+2}^{(2y_{k+2}-1)b}(q_{k+2}-1)^{b}} = \dfrac{q_{k+2}^{\alpha_{k+2}}}{A_{2}},
	\end{equation*}
	where $A_{2} = (q_{k+2}-1)^{b-1}$. If $2 \nmid \alpha_{k+2}$, then $x_{k+2} = (\alpha_{k+2} + 1)/2$ is a positive integer. Let $m_{k+2} = q_{k+2}^{x_{k+2}}$ and $n_{k+2} = 1$. Then we have
	\begin{equation*}
		\dfrac{\varphi\left(m_{k+2}^{2}\right)}{\left(\varphi\left(n_{k+2}^{2}\right)\right)^{b}} = \dfrac{\varphi\left(q_{k+2}^{2x_{k+2}}\right)}{\left(\varphi\left(1^{2}\right)\right)^{b}} = q_{k+2}^{2x_{k+2}-1}(q_{k+2} - 1) = \dfrac{q_{k+2}^{\alpha_{k+2}}}{A_{2}},
	\end{equation*}
	where $A_{2} = 1/(q_{k+2}-1)$.
	
	\textbf{Step $k + 3$:} Let $q_{k+3} = p_{s-2}$ and $\alpha_{k+3} = v_{q_{k+3}}(qA_{0}A_{1}A_{2})$. Since
	\begin{equation*}
		\begin{aligned}
			v_{q_{k+3}}(A)\geqslant t(b-1) & > \max_{1\leqslant j \leqslant s}\left\{v_{p_{j}}\left(\dfrac{1}{q}\prod_{i=1}^{s}(p_{i}-1)\right)\right\} \\
			& > v_{q_{k+3}}\left((q_{k+1}-1)(q_{k+2}-1)\right)-v_{q_{k+3}}(q),
		\end{aligned}
	\end{equation*}
	it follows that
	\begin{equation*}
		\alpha_{k+3} = v_{q_{k+3}}(qA_{0}A_{1}A_{2}) \geqslant v_{q_{k+3}}(q) + v_{q_{k+3}}(A) - v_{q_{k+3}}\left((q_{k+1} - 1)(q_{k+2}-1)\right) > 0.
	\end{equation*}
	Note that $b$ is an odd integer. If $2 \mid \alpha_{k+3}$, then there are positive integers $x_{k+3}$ and $y_{k+3}$ such that
	\begin{equation*}
		(2x_{k+3}-1) - (2y_{k+3}-1)b = \alpha_{k+3}.
	\end{equation*}
	Let $m_{k+3} = q_{k+3}^{x_{k+3}}$ and $n_{k+3} = q_{k+3}^{y_{k+3}}$. Then we have
	\begin{equation*}
		\dfrac{\varphi\left(m_{k+3}^{2}\right)}{\left(\varphi\left(n_{k+3}^{2}\right)\right)^{b}} = \dfrac{\varphi\left(q_{k+3}^{2x_{k+3}}\right)}{\left(\varphi\left(q_{k+3}^{2y_{k+3}}\right)\right)^{b}} = \dfrac{q_{k+3}^{2x_{k+3}-1}(q_{k+3} - 1)}{q_{k+3}^{(2y_{k+3}-1)b}(q_{k+3}-1)^{b}} = \dfrac{q_{k+3}^{\alpha_{k+3}}}{A_{3}},
	\end{equation*}
	where $A_{3} = (q_{k+3}-1)^{b-1}$. If $2 \nmid \alpha_{k+3}$, then $x_{k+3} = (\alpha_{k+3} + 1)/2$ is a positive integer. Let $m_{k+3} = q_{k+3}^{x_{k+3}}$ and $n_{k+3} = 1$. Then we have
	\begin{equation*}
		\dfrac{\varphi\left(m_{k+3}^{2}\right)}{\left(\varphi\left(n_{k+3}^{2}\right)\right)^{b}} = \dfrac{\varphi\left(q_{k+3}^{2x_{k+3}}\right)}{\left(\varphi\left(1^{2}\right)\right)^{b}} = q_{k+3}^{2x_{k+3}-1}(q_{k+3} - 1) = \dfrac{q_{k+3}^{\alpha_{k+3}}}{A_{3}},
	\end{equation*}
	where $A_{3} = 1/(q_{k+3}-1)$.
	
	Continuing this procedure, we may end at step $k + s$.
	
	\textbf{Step $k + s$:} Let $q_{k+s} = p_{1} = 2$ and $\alpha_{k+s} = v_{2}\left(q\prod_{i=0}^{s-1}A_{i}\right)$. Since
	\begin{equation*}
		v_{2}(A_{0})\geqslant t(b-1) > \max_{1\leqslant j \leqslant s}\left\{v_{p_{j}}\left(\dfrac{1}{q}\prod_{i=1}^{s}(p_{i}-1)\right)\right\} \geqslant v_{2}\left(\prod_{i=1}^{s-1}(q_{k+i} - 1)\right)-v_{2}(q),
	\end{equation*}
	it follows that
	\begin{equation*}
		\alpha_{k+s} = v_{2}\left(q\prod_{i=0}^{s-1}A_{i}\right) \geqslant v_{2}(q) + v_{2}(A_{0}) - v_{2}\left(\prod_{i=1}^{s-1}(q_{k+i} - 1)\right) > 0.
	\end{equation*}
	Note that $b$ is an odd integer. If $2 \mid \alpha_{k+s}$, then there are positive integers $x_{k+s}$ and $y_{k+s}$ such that
	\begin{equation*}
		(2x_{k+s}-1) - (2y_{k+s}-1)b = \alpha_{k+s}.
	\end{equation*}
	Let $m_{k+s} = 2^{x_{k+s}}$ and $n_{k+s} = 2^{y_{k+s}}$. Then we have
	\begin{equation*}
		\dfrac{\varphi\left(m_{k+s}^{2}\right)}{\left(\varphi\left(n_{k+s}^{2}\right)\right)^{b}} = \dfrac{\varphi\left(2^{2x_{k+s}}\right)}{\left(\varphi\left(2^{2y_{k+s}}\right)\right)^{b}} = \dfrac{2^{2x_{k+s}-1}(2 - 1)}{2^{(2y_{k+s}-1)b}(2-1)^{b}} = 2^{\alpha_{k+s}}.
	\end{equation*}
	If $2 \nmid \alpha_{k+s}$, then $x_{k+s} = (\alpha_{k+s} + 1)/2$ is a positive integer. Let $m_{k+s} = 2^{x_{k+s}}$ and $n_{k+s} = 1$. Then we have
	\begin{equation*}
		\dfrac{\varphi\left(m_{k+s}^{2}\right)}{\left(\varphi\left(n_{k+s}^{2}\right)\right)^{b}} = \dfrac{\varphi\left(2^{2x_{k+s}}\right)}{\left(\varphi\left(1^{2}\right)\right)^{b}} = 2^{2x_{k+s}-1}(2 - 1) = 2^{\alpha_{k+s}}.
	\end{equation*}	
	
	Let $m = \prod_{i=1}^{k+s}m_{i}$ and $n = \prod_{i=1}^{k+s}n_{i}$. Then we have
	\begin{equation*}
		\dfrac{\varphi(m^{2})}{\left(\varphi(n^{2})\right)^{b}} = \prod_{i=1}^{k+s}\dfrac{\varphi\left(m_{i}^{2}\right)}{\left(\varphi\left(n_{i}^{2}\right)\right)^{b}} = \dfrac{1}{A}\cdot\prod_{i=k+1}^{k+s}\dfrac{\varphi\left(m_{i}^{2}\right)}{\left(\varphi\left(n_{i}^{2}\right)\right)^{b}} = \dfrac{\prod_{i=1}^{s}p_{k+i}^{\alpha_{k+i}}}{\prod_{i=0}^{s-1}A_{i}} = q.
	\end{equation*}
	This completes the proof.

\emph{Proof of Theorem \ref{Thm2}.} Let $q = u/v$ be a positive rational number, where $u, v\in\mathbb{N}$ and $\gcd(u, v) = 1$. Suppose that
	\begin{equation}
		kluv = \prod_{i=1}^{s}p_{i}^{\alpha_{i}}, ~\text{where}~ p_{1}<p_{2}<\dots<p_{s} ~\text{are primes and}~ \alpha_{i}\in\mathbb{N},
	\end{equation}
	be the standard factorization of $kluv$. Let $d = \prod_{i=1}^{s}p_{i}^{\delta_{i}}$, where
	\begin{equation*}
		\delta_{i} = \left\{
		\begin{array}{ll}
			1, & ~\text{if}~ 2\nmid \alpha_{i}, \\
			0, & ~\text{if}~ 2\mid \alpha_{i}.
		\end{array}
		\right.
	\end{equation*}
	
	We distinguish two cases as follows.
	
	\textbf{Case 1:} $d = 1$
	
	In this case, $\delta_{i} = 0$, that is, $2 \mid \alpha_{i}$ for every $i = 1, 2, \dots, s$. This implies that $kluv = c^{2}$ for some $c\in\mathbb{N}$. Thus, we have $q = u/v = c^{2}/(klv^{2})$. By Dirichlet's theorem on primes in an arithmetic progression, there is a positive integer $t$ such that $1 + v^{2}c^{2}lt$ is a prime greater than $\max\{2, k\}$. Let $m = 2v^{2}c^{2}lt + 1$ and $n = 2v^{3}cltk$. Note that $\rad(4v^{2}c^{2}ltk) = \rad(4v^{6}c^{2}l^{3}t^{2}k^{2})$. Then we have
	\begin{equation*}
		\begin{aligned}
			\dfrac{\varphi(k(m^{2}-1))}{\varphi(ln^{2})}
			& = \dfrac{\varphi(4v^{2}c^{2}ltk(1+v^{2}c^{2}lt))}{\varphi(4v^{6}c^{2}l^{3}t^{2}k^{2})} \\
			& = \dfrac{\varphi(4v^{2}c^{2}ltk)}{\varphi(4v^{6}c^{2}l^{3}t^{2}k^{2})}\cdot v^{2}c^{2}lt = \dfrac{4v^{2}c^{2}ltk}{4v^{6}c^{2}l^{3}t^{2}k^{2}}\cdot v^{2}c^{2}lt = \dfrac{c^{2}}{klv^{2}} = q.
		\end{aligned}
	\end{equation*}
	
	\textbf{Case 2:} $d > 1$
	
	In this case, we have $\alpha_{i}\equiv\delta_{i} ~(\Mod 2)$, that is, $v_{p_{i}}(kluv) \equiv v_{p_{i}}(d) ~(\Mod 2)$ for every $i = 1, 2, \dots, s$. This implies that
	\begin{equation*}
		v_{p_{i}}(k) + v_{p_{i}}(d) - v_{p_{i}}(l) \equiv v_{p_{i}}(u) - v_{p_{i}}(v) ~(\Mod 2).
	\end{equation*}
	Hence, for every $i = 1, 2, \dots, s$, there are positive integers $x_{i}$ and $y_{i}$ such that
	\begin{equation*}
		v_{p_{i}}(k) + v_{p_{i}}(d) + 2x_{i} - v_{p_{i}}(l) - 2y_{i} = v_{p_{i}}(u) - v_{p_{i}}(v).
	\end{equation*}
	Let $M = \prod_{i=1}^{s}p_{i}^{x_{i}}$ and $N = \prod_{i=1}^{s}p_{i}^{y_{i}}$. Since $d$ is square-free, it follows that the Pell's equation $x^{2} - dM^{2}y^{2} = 1$ has infinitely many solutions in positive integers $x$ and $y$. Then we have $m_{0}^{2} - dM^{2}n_{0} = 1$ for some $m_{0}, n_{0}\in\mathbb{N}$. Let $m = m_{0}$ and $n = n_{0}N$. Since $\rad(kdM^{2}n_{0}^{2}) = \rad(lN^{2}n_{0}^{2})$, we obtain
	\begin{equation*}
		\dfrac{\varphi(k(m^{2}-1))}{\varphi(ln^{2})} = \dfrac{\varphi(kdM^{2}n_{0}^{2})}{\varphi(lN^{2}n_{0}^{2})} = \dfrac{kdM^{2}n_{0}^{2}}{lN^{2}n_{0}^{2}} = \dfrac{kdM^{2}}{lN^{2}}.
	\end{equation*}
	Therefore, for every $i = 1, 2,\dots, s$, we have
	\begin{equation*}
		\begin{aligned}
			v_{p_{i}}\left(\dfrac{\varphi(k(m^{2}-1))}{\varphi(ln^{2})} \right)
			& = v_{p_{i}}\left(\dfrac{kdM^{2}}{lN^{2}}\right) = v_{p_{i}}(k) + v_{p_{i}}(d) + 2x_{i} - v_{p_{i}}(l) - 2y_{i} \\
			& = v_{p_{i}}(u) - v_{p_{i}}(v) = v_{p_{i}}\left(\dfrac{u}{v}\right) = v_{p_{i}}(q).
		\end{aligned}
	\end{equation*}
	It follows that $q = \varphi(k(m^{2}-1))/\varphi(ln^{2})$. This completes the proof.

\begin{example}
	Find a positive integer pair $(m, n)$ such that
	\begin{itemize}
		\item[\rm (i)] $\dfrac{5}{12} = \dfrac{\varphi(m^{2})}{\left(\varphi(n^{2})\right)^{3}}$;
		\item[\rm (ii)] $\dfrac{5}{12} = \dfrac{\varphi(m^{2})}{\left(\varphi(n^{2})\right)^{5}}$.
	\end{itemize}
\end{example}

\emph{Solution.} (i) Note that $12 \cdot(2-1)\cdot(3-1)\cdot(5-1)/5 = 2^{5}\cdot 3^{1}\cdot 5^{-1}$ and
	\begin{equation*}
		t\cdot (3 - 1) > \max\{5, 1, -1\} = 5.
	\end{equation*}
	It follows that $t > 5/2$. We may let $t = 3$. Since $3456001$ is a prime such that $3456001\equiv 1 (\Mod (5!)^{3})$, we have
	\begin{equation*}
		\begin{aligned}
			\dfrac{5}{12} & = \dfrac{\varphi(3456001^{4})}{\left(\varphi(3456001^{2})\right)^{3}} \cdot (3456001 - 1)^{2}\cdot\dfrac{5}{12} = \dfrac{\varphi(3456001^{4})}{\left(\varphi(3456001^{2})\right)^{3}}\cdot 5^{7}\cdot 3^{5}\cdot 2^{18} \\
			& = \dfrac{\varphi(3456001^{4})}{\left(\varphi(3456001^{2})\right)^{3}}\cdot \dfrac{\varphi(5^{8})}{\left(\varphi(1^{2})\right)^{3}}\cdot \dfrac{3^{5}\cdot 2^{18}}{5 - 1} = \dfrac{\varphi(3456001^{4})}{\left(\varphi(3456001^{2})\right)^{3}}\cdot \dfrac{\varphi(5^{8})}{\left(\varphi(1^{2})\right)^{3}}\cdot 3^{5}\cdot 2^{16} \\
			& =  \dfrac{\varphi(3456001^{4})}{\left(\varphi(3456001^{2})\right)^{3}}\cdot \dfrac{\varphi(5^{8})}{\left(\varphi(1^{2})\right)^{3}}\cdot \dfrac{\varphi(3^{6})}{\left(\varphi(1^{2})\right)^{3}}\cdot 2^{15} \\
			& = \dfrac{\varphi(3456001^{4})}{\left(\varphi(3456001^{2})\right)^{3}}\cdot \dfrac{\varphi(5^{8})}{\left(\varphi(1^{2})\right)^{3}}\cdot \dfrac{\varphi(3^{6})}{\left(\varphi(1^{2})\right)^{3}}\cdot \dfrac{\varphi(2^{16})}{\left(\varphi(1^{2})\right)^{3}} \\
			& = \dfrac{\varphi\left((3456001^{2}\cdot 5^{4}\cdot 3^{3}\cdot 2^{8})^{2}\right)}{\left(\varphi\left(3456001^{2}\right)\right)^{3}}.
		\end{aligned}
	\end{equation*}
	(ii) Note that $12 \cdot(2-1)\cdot(3-1)\cdot(5-1)/5 = 2^{5}\cdot 3^{1}\cdot 5^{-1}$ and
	\begin{equation*}
		t\cdot (5 - 1) > \max\{5, 1, -1\} = 5.
	\end{equation*}
	It follows that $t > 5/4$. We may let $t = 2$. Since $14401$ is a prime such that $14401\equiv 1 (\Mod (5!)^{2})$, we have
	\begin{equation*}
		\begin{aligned}
			\dfrac{5}{12} & = \dfrac{\varphi(14401^{6})}{\left(\varphi(14401^{2})\right)^{5}} \cdot (14401 - 1)^{4}\cdot\dfrac{5}{12} = \dfrac{\varphi(14401^{6})}{\left(\varphi(14401^{2})\right)^{5}}\cdot 5^{9}\cdot 3^{7}\cdot 2^{22} \\
			& = \dfrac{\varphi(14401^{6})}{\left(\varphi(14401^{2})\right)^{5}}\cdot \dfrac{\varphi(5^{10})}{\left(\varphi(1^{2})\right)^{5}}\cdot \dfrac{3^{7}\cdot 2^{22}}{5 - 1} = \dfrac{\varphi(14401^{6})}{\left(\varphi(14401^{2})\right)^{5}}\cdot \dfrac{\varphi(5^{10})}{\left(\varphi(1^{2})\right)^{5}}\cdot 3^{7}\cdot 2^{20} \\
			& =  \dfrac{\varphi(14401^{6})}{\left(\varphi(14401^{2})\right)^{5}}\cdot \dfrac{\varphi(5^{10})}{\left(\varphi(1^{2})\right)^{5}}\cdot \dfrac{\varphi(3^{8})}{\left(\varphi(1^{2})\right)^{5}}\cdot 2^{19} \\
			& = \dfrac{\varphi(14401^{6})}{\left(\varphi(14401^{2})\right)^{5}}\cdot \dfrac{\varphi(5^{10})}{\left(\varphi(1^{2})\right)^{5}}\cdot \dfrac{\varphi(3^{8})}{\left(\varphi(1^{2})\right)^{5}}\cdot \dfrac{\varphi(2^{20})}{\left(\varphi(1^{2})\right)^{5}} \\
			& = \dfrac{\varphi\left((14401^{3}\cdot 5^{5}\cdot 3^{4}\cdot 2^{10})^{2}\right)}{\left(\varphi\left(14401^{2}\right)\right)^{5}}.\hspace{6.4cm}\qedhere
		\end{aligned}
	\end{equation*}

\begin{remark}
\rm
Using the method provided by the proof of Theorem \ref{Thm1}, we must find the positive integer pair $(m, n)$ that satisfies the condition. However, such pairs are not necessarily optimal. For example, we have
	\begin{equation*}
		\begin{aligned}
			\dfrac{5}{12} & = \dfrac{\varphi(241^{4})}{\left(\varphi(241^{2})\right)^{3}} \cdot (241 - 1)^{2}\cdot\dfrac{5}{12} = \dfrac{\varphi(241^{4})}{\left(\varphi(241^{2})\right)^{3}}\cdot 5^{3}\cdot 3^{1}\cdot 2^{6} \\
			& = \dfrac{\varphi(241^{4})}{\left(\varphi(241^{2})\right)^{3}}\cdot \dfrac{\varphi(5^{4})}{\left(\varphi(1^{2})\right)^{3}}\cdot \dfrac{3^{1}\cdot 2^{6}}{5 - 1} = \dfrac{\varphi(241^{4})}{\left(\varphi(241^{2})\right)^{3}}\cdot \dfrac{\varphi(5^{4})}{\left(\varphi(1^{2})\right)^{3}}\cdot 3^{1}\cdot 2^{4} \\
			& =  \dfrac{\varphi(241^{4})}{\left(\varphi(241^{2})\right)^{3}}\cdot \dfrac{\varphi(5^{4})}{\left(\varphi(1^{2})\right)^{3}}\cdot \dfrac{\varphi(3^{2})}{\left(\varphi(1^{2})\right)^{3}}\cdot 2^{3} \\
			& = \dfrac{\varphi(241^{4})}{\left(\varphi(241^{2})\right)^{3}}\cdot \dfrac{\varphi(5^{4})}{\left(\varphi(1^{2})\right)^{3}}\cdot \dfrac{\varphi(3^{2})}{\left(\varphi(1^{2})\right)^{3}}\cdot \dfrac{\varphi(2^{4})}{\left(\varphi(1^{2})\right)^{3}} \\
			& = \dfrac{\varphi\left((241^{2}\cdot 5^{2}\cdot 3^{1}\cdot 2^{2})^{2}\right)}{\left(\varphi\left(241^{2}\right)\right)^{3}}.
		\end{aligned}
	\end{equation*}
	and
	\begin{equation*}
		\begin{aligned}
			\dfrac{5}{12} & = \dfrac{\varphi(61^{6})}{\left(\varphi(61^{2})\right)^{5}} \cdot (61 - 1)^{4}\cdot\dfrac{5}{12} = \dfrac{\varphi(61^{6})}{\left(\varphi(61^{2})\right)^{5}}\cdot 5^{5}\cdot 3^{3}\cdot 2^{6} \\
			& = \dfrac{\varphi(61^{6})}{\left(\varphi(61^{2})\right)^{5}}\cdot \dfrac{\varphi(5^{6})}{\left(\varphi(1^{2})\right)^{5}}\cdot \dfrac{3^{3}\cdot 2^{6}}{5 - 1} = \dfrac{\varphi(61^{6})}{\left(\varphi(61^{2})\right)^{5}}\cdot \dfrac{\varphi(5^{6})}{\left(\varphi(1^{2})\right)^{5}}\cdot 3^{3}\cdot 2^{4} \\
			& =  \dfrac{\varphi(61^{6})}{\left(\varphi(61^{2})\right)^{5}}\cdot \dfrac{\varphi(5^{6})}{\left(\varphi(1^{2})\right)^{5}}\cdot \dfrac{\varphi(3^{4})}{\left(\varphi(1^{2})\right)^{5}}\cdot 2^{3} \\
			& = \dfrac{\varphi(61^{6})}{\left(\varphi(61^{2})\right)^{5}}\cdot \dfrac{\varphi(5^{6})}{\left(\varphi(1^{2})\right)^{5}}\cdot \dfrac{\varphi(3^{4})}{\left(\varphi(1^{2})\right)^{5}}\cdot \dfrac{\varphi(2^{4})}{\left(\varphi(1^{2})\right)^{5}} \\
			& = \dfrac{\varphi\left((61^{3}\cdot 5^{3}\cdot 3^{2}\cdot 2^{2})^{2}\right)}{\left(\varphi\left(61^{2}\right)\right)^{5}}.
		\end{aligned}
	\end{equation*}
	
\end{remark}

\begin{example}
	Find a positive integer pair $(m, n)$ such that
	\begin{itemize}
		\item[\rm (i)] $\dfrac{5}{12} = \dfrac{\varphi(15(m^{2} - 1))}{\varphi(2n^{2})}$;
		\item[\rm (ii)] $\dfrac{5}{12} = \dfrac{\varphi(3(m^{2} - 1))}{\varphi(5n^{2})}$.
	\end{itemize}
\end{example}

\emph{Solution.}
	(i) Note that $2^{5}\cdot 3^{4}\cdot 5^{2} \cdot 7^{2}\cdot 11^{2} = 19601^{2} - 1$. Then we have
	\begin{equation*}
		\begin{aligned}
			\dfrac{5}{12} & = \dfrac{5}{2^{2}\cdot 3} = \dfrac{2^{5}\cdot 3^{5}\cdot 5^{3} \cdot 7^{2}\cdot 11^{2}}{2^{7}\cdot 3^{6}\cdot 5^{2} \cdot 7^{2}\cdot 11^{2}} = \dfrac{\varphi(2^{5}\cdot 3^{5}\cdot 5^{3} \cdot 7^{2}\cdot 11^{2})}{\varphi(2^{7}\cdot 3^{6}\cdot 5^{2} \cdot 7^{2}\cdot 11^{2})} \\
			& = \dfrac{\varphi(15\cdot(2^{5}\cdot 3^{4}\cdot 5^{2} \cdot 7^{2}\cdot 11^{2}))}{\varphi(2\cdot(2^{3}\cdot 3^{3}\cdot 5 \cdot 7\cdot 11)^{2})} = \dfrac{\varphi(15\cdot(19601^{2}-1))}{\varphi(2\cdot 83160^{2})}.
		\end{aligned}
	\end{equation*}
	(ii) Since $5185001 = 2^{9}\cdot 3^{4}\cdot 5^{3} + 1$, we obtain
	\begin{equation*}
		(2^{10}\cdot 3^{4}\cdot 5^{3} + 1)^{2} - 1 = 2^{11} \cdot 3^{4} \cdot 5^{3} \cdot (2^{9}\cdot 3^{4}\cdot 5^{3} + 1) = 2^{11} \cdot 3^{4} \cdot 5^{3} \cdot 5184001.
	\end{equation*}
	Note that $5184001$ is a prime. Then
	\begin{equation*}
		\begin{aligned}
			\dfrac{5}{12}
			& = \dfrac{5}{2^{2}\cdot 3} = \dfrac{2^{20}\cdot 3^{9}\cdot 5^{6}}{2^{22}\cdot 3^{10}\cdot 5^{5}} = \dfrac{2^{11}\cdot 3^{5}\cdot 5^{3}}{2^{22}\cdot 3^{10}\cdot 5^{5}}\cdot 2^{9}\cdot 3^{4}\cdot 5^{3} \\
			& = \dfrac{\varphi(2^{11}\cdot 3^{5}\cdot 5^{3})}{\varphi(2^{22}\cdot 3^{10}\cdot 5^{5}}\cdot \varphi(5184001) = \dfrac{\varphi(2^{11}\cdot 3^{5}\cdot 5^{3}\cdot 5184001)}{\varphi(2^{22}\cdot 3^{10}\cdot 5^{5})} \\
			& = \dfrac{\varphi(3\cdot(2^{11}\cdot 3^{4}\cdot 5^{3}\cdot 5184001))}{\varphi(5\cdot(2^{22}\cdot 3^{10}\cdot 5^{4}))} = \dfrac{\varphi(3\cdot\left((2^{10}\cdot 3^{4}\cdot 5^{3} + 1)^{2} - 1\right))}{\varphi(5\cdot(2^{11}\cdot 3^{5}\cdot 5^{2})^{2})}.\hspace{0.9cm}\qedhere
		\end{aligned}
	\end{equation*}

\section{Further Researches} \label{Sec3}

In this section, we further discuss the problem in Section \ref{Sec1}. Let $\Gamma$ be the set of all positive integer quadruples $(a, b, r, s)$ such that
\begin{equation*}
	\mathbb{Q^{+}} = \left\{\left.\dfrac{(\varphi(m^{r}))^{a}}{(\varphi(n^{s}))^{b}}\right| m, n\in\mathbb{N} \right\}.
\end{equation*}
It is obvious that $\gcd(a, b) = 1$. Moreover, we have the following facts.

\begin{fact}[2022, Vu \cite{V2022}]
	If $\gcd(ar, bs) = 1$, then $(a, b, r, s) \in\Gamma$.
\end{fact}

\begin{fact} \label{F2}
	$(1, 1, r, s)\in\Gamma$ if and only if $r = s = 2$.
\end{fact}

\begin{prof}
	This fact follows by the work of Krachun and Sun \cite{KS2020} and Li, Yuan and Bai \cite{LYB2020}.
\end{prof}

\begin{fact} \label{F3}
	If $\gcd(ar, bs) = d > 1$ and $d \nmid (a - b)$, then $(a, b, r, s)\notin\Gamma$.
\end{fact}

\begin{prof}
	We have just seen that the statement holds for $\gcd(a, b) > 1$. Now we consider $\gcd(a, b) = 1$, which implies that at least one of the positive integers $r$ and $s$ is greater than $1$. Without loss of generality, we may assume that $r > 1$ when $d > 2$. When $d = 2$, we may assume that $2 \nmid a$ and $2 \mid b$, which implies that $r > 1$.
	
	When $d > 2$ or $d = 2$ with $2 \nmid a$ and $2 \mid b$, we can choose a positive integer $k$ such that $d \nmid (a + k)$ and $d \nmid (a - b + k)$. It suffices to show that there are no positive integers $m$ and $n$ such that
	\begin{equation} \label{eq1}
		2^{k} = \dfrac{(\varphi(m^{r}))^{a}}{(\varphi(n^{s}))^{b}}.
	\end{equation}
	
	Clearly, it is impossible when $m = n = 1$. Let $p$ and $q$ be the maximal prime factor of $m$ and $n$, respectively.
	
	If $m = 1$ and $n > 1$, then $v_{2}\left(\left(\varphi(m^{r})\right)^{a}/\left(\varphi(n^{s})\right)^{b}\right) = -v_{2}\left(\varphi(n^s)\right)b \leqslant 0$. It follows that
	\begin{equation*}
		v_{2}\left(\dfrac{(\varphi(m^{r}))^{a}}{(\varphi(n^{s}))^{b}}\right) = -v_{2}\left(\varphi(n^s)\right)b \neq v_{2}\left(2^{k}\right),
	\end{equation*}
	which contradicts with \eqref{eq1}.
	
	If $m > 1$ and $n = 1$, since $r > 1$ and $d \nmid (a + k)$, then there are no positive integer $x$ such that $(rx - 1)a = v_{p}(2^{k})$. It follows that
	\begin{equation*}
		v_{p}\left(\dfrac{(\varphi(m^{r}))^{a}}{(\varphi(n^{s}))^{b}}\right) = (rv_{p}(m) - 1)a \neq v_{p}\left(2^{k}\right),
	\end{equation*}
	which contradicts with \eqref{eq1}.
	
	If $m > 1$ and $n > 1$, then we distinguish three cases as follows.
	
	\textbf{Case 1: } $p > q \geqslant 2$
	
	In this case, since $r > 1$, we have
	\begin{equation*}
		v_{p}\left(\dfrac{(\varphi(m^{r}))^{a}}{(\varphi(n^{s}))^{b}}\right) = (rv_{p}(m)- 1)a \neq v_{p}(2^{k}),
	\end{equation*}
	which contradicts with \eqref{eq1}.
	
	\textbf{Case 2: } $q > p \geqslant 2$
	
	In this case, if $s > 1$, then we have
	\begin{equation*}
		v_{q}\left(\dfrac{(\varphi(m^{r}))^{a}}{(\varphi(n^{s}))^{b}}\right) = -(sv_{q}(n)- 1)b \neq v_{q}(2^{k}),
	\end{equation*}
	which contradicts with \eqref{eq1}.  If $s = 1$, then $d \mid b$ and $d \nmid a$. Since $d \nmid a$ and $d \nmid (a + k)$, there are no positive integers $x$ and $y$ such that
	\begin{equation*}
		(rx - 1)a - by = v_{p}(2^{k}).
	\end{equation*}
	It follows that
	\begin{equation*}
		v_{p}\left(\dfrac{(\varphi(m^{r}))^{a}}{(\varphi(n^{s}))^{b}}\right) = (rv_{p}(m) - 1)a -v_{p}\left(\varphi(n^{s})\right)b \neq v_{p}(2^{k}),
	\end{equation*}
	which contradicts with \eqref{eq1}.
	
	\textbf{Case 3: } $p = q \geqslant 2$
	
	In this case, since $d \nmid (a-b)$ and $d \nmid (a-b+k)$, then there are no positive integers $x$ and $y$ such that
	\begin{equation*}
		(rx - 1)a - (sy - 1)b = v_{p}(2^{k}).
	\end{equation*}
	Thus, we have
	\begin{equation*}
		v_{p}\left(\dfrac{(\varphi(m^{r}))^{a}}{(\varphi(n^{s}))^{b}}\right) = (rv_{p}(n) - 1)a -(sv_{p}(n)- 1)b \neq v_{p}(2^{k}).
	\end{equation*}
	which contradicts with \eqref{eq1} and the proof is completed.
\end{prof}

\begin{fact} \label{F4}
	If $\gcd(ar, bs) = d > 2 $ and $d \mid (a - b)$, then $(a, b, r, s)\notin\Gamma$.
\end{fact}

\begin{prof}
	We have just seen that the statement holds for $\gcd(a, b) > 1$. Now we consider $\gcd(a, b) = 1$. Since $d \mid (a-b)$, it follows that $d\nmid a$ and $d\nmid b$. Thus, $r > 1$ and $s > 1$. By Fact \ref{F2}, we have $a\neq b$. Without loss of generality, we may assume that $a < b$.
	
	Since $d > 2$, we can choose a positive integer $k$ such that $d \nmid k$ and $d \nmid (a + k)$. It suffices to show that there are no positive integers $m$ and $n$ such that
	\begin{equation} \label{eq2}
		2^{k} = \dfrac{(\varphi(m^{r}))^{a}}{(\varphi(n^{s}))^{b}}.
	\end{equation}
	
	Clearly, it is impossible when $m = n = 1$. Let $p$ and $q$ be the maximal prime factor of $m$ and $n$, respectively.
	
	If $m = 1$ and $n > 1$, then $v_{2}\left(\left(\varphi(m^{r})\right)^{a}/\left(\varphi(n^{s})\right)^{b}\right) = -v_{2}\left(\varphi(n^s)\right)b \leqslant 0$. It follows that
	\begin{equation*}
		v_{2}\left(\dfrac{(\varphi(m^{r}))^{a}}{(\varphi(n^{s}))^{b}}\right) = -v_{2}\left(\varphi(n^s)\right)b \neq v_{2}\left(2^{k}\right),
	\end{equation*}
	which contradicts with \eqref{eq2}.
	
	If $m > 1$ and $n = 1$, since $r > 1$ and $d \nmid (a + k)$, then there are no positive integer $x$ such that $(rx - 1)a = v_{p}(2^{k})$. It follows that
	\begin{equation*}
		v_{p}\left(\dfrac{(\varphi(m^{r}))^{a}}{(\varphi(n^{s}))^{b}}\right) = (rv_{p}(m) - 1)a \neq v_{p}\left(2^{k}\right),
	\end{equation*}
	which contradicts with \eqref{eq2}.
	
	If $m > 1$ and $n > 1$, then we distinguish three cases as follows.
	
	\textbf{Case 1: } $p > q \geqslant 2$
	
	In this case, since $r > 1$, we have
	\begin{equation*}
		v_{p}\left(\dfrac{(\varphi(m^{r}))^{a}}{(\varphi(n^{s}))^{b}}\right) = (rv_{p}(m)- 1)a \neq v_{p}(2^{k}),
	\end{equation*}
	which contradicts with \eqref{eq2}.
	
	\textbf{Case 2: } $q > p \geqslant 2$
	
	In this case, since $s > 1$, we have
	\begin{equation*}
		v_{q}\left(\dfrac{(\varphi(m^{r}))^{a}}{(\varphi(n^{s}))^{b}}\right) = -(sv_{q}(n)- 1)b \neq v_{q}(2^{k}),
	\end{equation*}
	which contradicts with \eqref{eq2}.
	
	\textbf{Case 3: } $p = q \geqslant 2$
	
	In this case, since $d \mid (a-b)$ and $d \nmid k$, it follows that $d \nmid a$, $d \nmid b$ and $d \nmid (a-b+k)$. If $p = 2$, then there are no positive integers $x$ and $y$ such that
	\begin{equation*}
		(rx-1)a - (sy-1)b = k.
	\end{equation*}
	It follows that
	\begin{equation*}
		v_{2}\left(\dfrac{(\varphi(m^{r}))^{a}}{(\varphi(n^{s}))^{b}}\right) = (rv_{2}(n) - 1)a -(sv_{2}(n)- 1)b \neq v_{2}(2^{k}),
	\end{equation*}
	which contradicts with \eqref{eq2}. If $p > 2$, then we may assume that $m = p_{1}^{x_{1}}m_{1}$ and $n = p_{1}^{y_{1}}n_{1}$, where $m_{1}, n_{1}\in\mathbb{N}$, $p_{1} = p$, $x_{1} = v_{p_{1}}(m)$ and $y_{1} = v_{p_{1}}(n)$. Thus, we have
	\begin{equation*}
		\dfrac{(\varphi(m^{r}))^{a}}{(\varphi(n^{s}))^{b}} = \dfrac{p_{1}^{\alpha_{1}}}{(p_{1} - 1)^{b-a}}\cdot\dfrac{(\varphi(m_{1}^{r}))^{a}}{(\varphi(n_{1}^{s}))^{b}},
	\end{equation*}
	where $\alpha_{1} = (rx_{1}-1)a - (sy_{1}-1)b$. By \eqref{eq2}, we have $\alpha_{1} = 0$. Let $p_{2}$ and $q_{2}$ be the maximal prime factor of $m_{1}$ and $n_{1}$, respectively. Since $d \nmid a$ and $d \nmid b$, there are no positive integers $x$ and $y$ such that
	\begin{equation*}
		(rx - 1)a - (b-a)v_{p_{2}}(p_{1}-1) = 0
	\end{equation*}
	and
	\begin{equation*}
		-(sy - 1)b - (b-a)v_{q_{2}}(p_{1}-1) = 0.
	\end{equation*}
	It follows that $p_{2} = q_{2}$.
	
	Let $m_{1} = p_{2}^{x_{2}}m_{2}$ and $n_{1} = p_{2}^{y_{2}}n_{2}$, where $m_{2}, n_{2}\in\mathbb{N}$, $x_{2} = v_{p_{2}}(m_{1})$ and $y_{2} = v_{p_{2}}(n_{1})$. Then we have
	\begin{equation*}
		\begin{aligned}
			\dfrac{(\varphi(m^{r}))^{a}}{(\varphi(n^{s}))^{b}} & = \dfrac{p_{1}^{\alpha_{1}}}{(p_{1} - 1)^{b-a}}\cdot\dfrac{(\varphi(m_{1}^{r}))^{a}}{(\varphi(n_{1}^{s}))^{b}}\\
			& = \dfrac{p_{1}^{\alpha_{1}}}{(p_{1} - 1)^{b-a}}\cdot\dfrac{p_{2}^{\alpha_{2}}}{(p_{2} - 1)^{b-a}}\cdot\dfrac{(\varphi(m_{2}^{r}))^{a}}{(\varphi(n_{2}^{s}))^{b}},
		\end{aligned}
	\end{equation*}
	where $\alpha_{2} = (rx_{2}-1)a - (sy_{2}-1)b$. By \eqref{eq2}, we have $\alpha_{2} = (b-a)v_{p_{2}}(p_{1} - 1)$. Let $p_{3}$ and $q_{3}$ be the maximal prime factor of $m_{2}$ and $n_{2}$, respectively. Since $d \nmid a$ and $d \nmid b$, there are no positive integers $x$ and $y$ such that
	\begin{equation*}
		(rx - 1)a - (b-a)v_{p_{3}}\left((p_{1}-1)(p_{2}-1)\right) = 0
	\end{equation*}
	and
	\begin{equation*}
		-(sy - 1)b - (b-a)v_{q_{3}}\left((p_{1}-1)(p_{2}-1)\right) = 0.
	\end{equation*}
	It follows that $p_{3} = q_{3}$.
	
	Let $m_{2} = p_{3}^{x_{3}}m_{3}$ and $n_{2} = p_{3}^{y_{3}}n_{3}$, where $m_{3}, n_{3}\in\mathbb{N}$, $x_{3} = v_{p_{3}}(m_{2})$ and $y_{3} = v_{p_{3}}(n_{2})$. Then we have
	\begin{equation*}
		\begin{aligned}
			\dfrac{(\varphi(m^{r}))^{a}}{(\varphi(n^{s}))^{b}} & = \dfrac{p_{1}^{\alpha_{1}}}{(p_{1} - 1)^{b-a}}\cdot\dfrac{(\varphi(m_{1}^{r}))^{a}}{(\varphi(n_{1}^{s}))^{b}}\\
			& = \dfrac{p_{1}^{\alpha_{1}}}{(p_{1} - 1)^{b-a}}\cdot\dfrac{p_{2}^{\alpha_{2}}}{(p_{2} - 1)^{b-a}}\cdot\dfrac{(\varphi(m_{2}^{r}))^{a}}{(\varphi(n_{2}^{s}))^{b}}\\
			& = \dfrac{p_{1}^{\alpha_{1}}}{(p_{1} - 1)^{b-a}}\cdot\dfrac{p_{2}^{\alpha_{2}}}{(p_{2} - 1)^{b-a}}\cdot\dfrac{p_{3}^{\alpha_{3}}}{(p_{3} - 1)^{b-a}}\cdot\dfrac{(\varphi(m_{3}^{r}))^{a}}{(\varphi(n_{3}^{s}))^{b}},
		\end{aligned}
	\end{equation*}
	where $\alpha_{3} = (rx_{3}-1)a - (sy_{3}-1)b$. By \eqref{eq2}, we have $\alpha_{3} = (b-a)v_{p_{3}}\left((p_{1} - 1)(p_{2}-1)\right)$.
	
	Continuing this procedure, we assume that $s$ is the maximal positive integer such that $p_{s} > 2$. Then we have
	\begin{equation*}
		m = 2^{\alpha}\prod_{i=1}^{s}p_{i}^{x_{i}} ~~\text{and}~~ n = 2^{\beta}\prod_{i=1}^{s}p_{i}^{y_{i}},
	\end{equation*}
	where $p = p_{1} > p_{2} > \dots > p_{s} > 2$ are primes and $\alpha$, $\beta$ are nonnegative integers. Thus, we have
	\begin{equation*}
		\dfrac{(\varphi(m^{r}))^{a}}{(\varphi(n^{s}))^{b}} = \dfrac{(\varphi(2^{r\alpha}))^{a}}{(\varphi(2^{s\beta}))^{b}}\cdot\prod_{i=1}^{s}\dfrac{p_{i}^{\alpha_{i}}}{(p_{i} - 1)^{b-a}} = \dfrac{(\varphi(2^{r\alpha}))^{a}}{(\varphi(2^{s\beta}))^{b}} \cdot \dfrac{1}{2^{t(b-a)}}, ~~\text{where}~~ t\in\mathbb{N}.
	\end{equation*}
	
	If $\alpha = 0$ and $\beta = 0$, since $d \nmid k$ or $a < b$, then we have
	\begin{equation*}
		v_{2}\left(\dfrac{(\varphi(m^{r}))^{a}}{(\varphi(n^{s}))^{b}}\right) = -t(b-a) \neq v_{2}(2^{k}),
	\end{equation*}
	which contradicts with \eqref{eq2}.
	
	If $\alpha > 0$ and $\beta = 0$, since $d \nmid (a + k)$, then we have
	\begin{equation*}
		v_{2}\left(\dfrac{(\varphi(m^{r}))^{a}}{(\varphi(n^{s}))^{b}}\right) = (r\alpha - 1)a -t(b-a) \neq v_{2}(2^{k}),
	\end{equation*}
	which contradicts with \eqref{eq2}.
	
	If $\alpha = 0$ and $\beta > 0$, since $s > 1$ and $a < b$, then we have
	\begin{equation*}
		v_{2}\left(\dfrac{(\varphi(m^{r}))^{a}}{(\varphi(n^{s}))^{b}}\right) = -(s\beta - 1)b -t(b-a) \neq v_{2}(2^{k}),
	\end{equation*}
	which contradicts with \eqref{eq2}.
	
	If $\alpha > 0$ and $\beta > 0$, since $d \nmid (a - b + k)$, then we have
	\begin{equation*}
		v_{2}\left(\dfrac{(\varphi(m^{r}))^{a}}{(\varphi(n^{s}))^{b}}\right) = (r\alpha - 1)a - (s\beta - 1)b -t(b-a) \neq v_{2}(2^{k}),
	\end{equation*}
	which contradicts with \eqref{eq2} and the proof is completed.
\end{prof}

The remaining case is $\gcd(ar, bs) = 2$ with $2\nmid ab$. In this case, Theorem \ref{Thm1} shows that $(a, 1, 2, 2), (1, b, 2, 2)\in\Gamma$ when $a$ and $b$ are odd integers greater than $1$. In fact, using the same methods provided by the proof of Theorem \ref{Thm1}, we can verify the following fact.

\begin{fact} \label{F5}
	$(a, 1, r, 2), (1, b, 2, s)\in\Gamma$, where $a, b$ are two odd positive integers greater than $1$ and $r, s$ are two even positive integers.
\end{fact}

Up to now, besides the quadruples listed in Fact \ref{F2} and Fact \ref{F5}, we still do not know whether there are other quadruples $(a, b, r, s)$ with $\gcd(ar, bs) = 2$ in $\Gamma$.  This now leads us to propose the following question.

\begin{question}
	Can every positive rational number $q$ can be written in the form
	\begin{equation*}
		q = \dfrac{\left(\varphi(m^{2})\right)^{3}}{\left(\varphi(n^{2})\right)^{5}}, ~\text{where}~ m, n\in\mathbb{N}~ ?
	\end{equation*}
\end{question}

Inspired by Theorem \ref{Thm2}, we propose the following open problem for further research.

\begin{question}
	Let $t$ be an integer which is not a square. Are there positive integer pairs $(k, l)$ such that every positive rational number $q$ can be written in the form
	\begin{equation*}
		q = \dfrac{\varphi(k(m^{2}-t))}{\varphi(ln^{2})}, ~\text{where}~ m, n\in\mathbb{N}~ ?
	\end{equation*}
\end{question}

For example, let $t = -1$, $k = 1$ and $l = 1$. Can every positive rational number $q$ be written in the form
\begin{equation*}
	q = \dfrac{\varphi(m^{2} + 1)}{\varphi(n^{2})}, ~\text{where}~ m, n\in\mathbb{N}~ ?
\end{equation*}
More specifically, for every positive integer $w$, can $2^{w}$ be written in the form
\begin{equation*}
	2^{w} = \dfrac{\varphi(m^{2} + 1)}{\varphi(n^{2})}, ~\text{where}~ m, n\in\mathbb{N}~ ?
\end{equation*}

\end{document}